\documentclass[11pt]{article}

\usepackage[utf8]{inputenc}
\usepackage[english]{babel}
\usepackage{amsfonts}
\usepackage{url}
\usepackage{pdfpages}
\usepackage{authblk}

\usepackage{amsmath}
\usepackage{amssymb}
\usepackage{mathtools}
\usepackage{amsthm}

\usepackage[outline]{contour} 

\usepackage[left=3.54cm,right=3.54cm,top=4cm,bottom=4cm]{geometry}

\usepackage{booktabs}  
\usepackage{graphicx}

\usepackage{tikz}
\usetikzlibrary{arrows,shapes,positioning}

\usepackage{enumitem}

\usepackage{fancyhdr}

\DeclareMathOperator*{\argmin}{argmin}
\DeclareMathOperator*{\argmax}{argmax}

\newtheorem{definition}{Definition}

\usepackage{hyperref}
\usepackage{hypcap}

\usepackage{algorithm}
\usepackage{algorithmic}
\usepackage[algo2e,ruled,vlined]{algorithm2e}

\usepackage{glossaries}

\usepackage{float}

\usepackage{xcolor}

\usepackage{multirow}

\def\xx{\boldsymbol{x}}

\def\ss{\boldsymbol{s}}
\def\yy{\boldsymbol{y}}

\def\mm{\boldsymbol{m}}

\title {A machine learning framework for neighbor generation in metaheuristic search} 

\author[1]{Defeng Liu, }
\author[1]{Vincent Perreault}
\author[1]{Alain Hertz}
\author[1,2]{Andrea Lodi}

\affil[1]{\small MAGI, Polytechnique Montreal, Montreal, Canada}
\affil[2]{\small Jacobs Technion-Cornell Institute, Cornell Tech and Technion - IIT, New York,  USA}

\date{\vspace{-5ex}}

\providecommand{\keywords}[1]{\textbf{\textit{Keywords:}} #1}

\begin{document}
\maketitle

\begin{abstract}

This paper presents a methodology for integrating machine learning techniques into metaheuristics for solving combinatorial optimization problems. Namely, we propose a general machine learning framework for neighbor generation in metaheuristic search. We first define an efficient neighborhood structure constructed by applying a transformation to a selected subset of variables from the current solution. Then, the key of the proposed methodology is to generate promising neighbors by selecting a proper subset of variables that contains a descent of the objective in the solution space. To learn a good \emph{variable selection} strategy, we formulate the problem as a \emph{classification} task that exploits structural information from the characteristics of the problem and from high-quality solutions. We validate our methodology on two metaheuristic applications: a \emph{Tabu Search} scheme for solving a Wireless Network Optimization problem and a \emph{Large Neighborhood Search} heuristic for solving Mixed-Integer Programs. The experimental results show that our approach is able to achieve a satisfactory trade-off between the exploration of a larger solution space and the exploitation of high-quality solution regions on both applications. 

\end{abstract}


\section{Introduction}
\label{sec: intro}

Combinatorial Optimization (CO) is an important class of optimization problems in Operations Research (OR) and Computer Science (CS). In general, a CO problem is defined by a set of decision variables, a constrained solution space and an objective function. The goal of CO is to find optimal solutions with respect to the objective in the solution space.

Classical methods for solving CO problems can be roughly divided into three classes: exact methods, heuristics and metaheuristics. Mixed-Integer Programming (MIP) is one of the main paradigms of exact methods for modeling complex CO problems. Over the last decades, there has been increasing interest in improving the ability to solve MIPs effectively. Modern MIP solvers incorporate a variety of complex algorithmic techniques, such as primal heuristics \cite{blum2003metaheuristics}, Branch and Bound (B\&B) \cite{land2010automatic}, cutting planes \cite{Bixby2007} and pre-processing, which results in complex and sophisticated software tools.

Exact methods are guaranteed to find optimal solutions as well as a proof of their optimality. On the one hand, due to the NP-hardness nature of many CO problems, solving them to optimality within an exact algorithm is still a very challenging task. On the other hand, in many practical applications, one is more interested in getting a good solution within a reasonable time rather than finding an optimal one. Those practical requirements have motivated the development of \emph{specific heuristics} and \emph{metaheuristics} (MHs). Specific heuristics are usually designed for solving a certain type of CO problem, whereas MHs are frameworks for designing heuristics for solving general CO problems, and provide guidelines to integrate basic heuristic concepts with high-level diversification strategies. 

It is worth noting that a lot of information can be produced and observed from MH processes, and therefore, a large volume of data can be collected. These data might provide valuable information about the optimization status of the process, the characteristics of the problem, the structures and properties of high-quality solutions in the solution regions being visited. However, such knowledge has not been fully exploited by traditional MH algorithms.

That can be viewed as a disadvantage when compared to application-targeted CO algorithms, in which exploiting the domain knowledge of the problem is generally preferable. Meanwhile, real-world CO problems have a rich structure. With similar instances repeatedly solved in many applications, statistical characteristics and patterns appear. This provides the opportunity for Machine Learning (ML) to extract structural properties of the problem from data and automatically produce learning-based MH strategies.

ML is a subfield of Artificial Intelligence (AI) that involves developing algorithms to learn knowledge from data and make predictions on new problems. In recent years, the application of ML techniques in CO became an emerging research area with quite a number of contributions for different purposes. On the one hand, some research has been devoted to develop heuristics for solving CO problems by ML, i.e., to perform ``end-to-end learning'' to directly generate good solutions for a CO problem. On the other hand, ML has been applied in combination with CO algorithms, where the learning based algorithms have the potential to achieve better performances, either because the current strategy of performing some auxiliary tasks is computationally expensive or because they are poorly understood from the mathematical viewpoint. For a detailed review of ``learn to optimize'',  the interested reader is referred to the survey \cite{bengio2021machine}.

Specifically for metaheuristics, ML techniques can be used to infer patterns from data generated from MH processes. Integrating the extracted knowledge into the search strategies can lead MHs to search the solution space more efficiently and significantly improve the current performance. Recently, the application of ML techniques for MHs has attracted increasing research interest and we refer the interested reader to the following surveys \cite{talbi2021machine, karimi2022machine}.

In this paper, we focus on the integration of ML techniques into MHs and propose a general learning-based framework for neighbor generation in MHs. We first define an efficient neighborhood structure by applying a transformation to a selected subset of variables from the current solution. Then, the problem is to determine how to select a subset of variables that leads to a promising neighborhood of solutions containing a descent of the objective in the solution space. By conducting a classification task, our method learns good variable selection policies from both structural characteristics of the problem and high-quality solutions. 
We will demonstrate the effectiveness of our approach on two applications: a Tabu Search scheme for solving a Wireless Network Optimization problem and a Large Neighborhood Search heuristic for solving MIPs.


\section{Background}
\label{sec:background}

In this section, we introduce the necessary background and notation.

\subsection{Combinatorial Optimization}
\label{sec:background_comh}

CO problems are a class of optimization problems with a set of decision variables and a defined solution space. Without loss of generality, a CO problem can be formulated into a constrained optimization problem as follows:
\begin{align}
 \min~~~ & f(\boldsymbol{x})\\
\text{s.t.}~~~  & g(\boldsymbol{x}) \le \boldsymbol{b},\\
& x_i \in \{0,1\}, ~\forall i\in \mathcal{B},\\
& x_j \in \mathbb{Z}^+, ~\forall j\in \mathcal{G}, \\
& x_k \ge 0, ~\forall k\in \mathcal{C},
\end{align}
where the index set $\mathcal{V}:=\{1,\ldots,n\}$ of decision variables is partitioned into $\mathcal{B}, \mathcal{G}, \mathcal{C}$, which are the index sets of binary, general integer and continuous variables, respectively.

Since many CO problems are NP-hard, determining optimal solutions by exact methods requires in the worst case exponential time and might be intractable, especially for large-size applications. In many practical applications where the CO problems are hard and complex, practitioners are often interested in finding good-quality solutions in an ``acceptable'' amount of computing time rather than solving the problem to optimality. Therefore, heuristic algorithms are developed to efficiently compute high-quality solutions. 

\subsection{Metaheuristics}
\label{sec: background_mh}
Metaheuristics are general framework strategies for designing heuristics for solving CO problems, and provide guidelines to integrate basic heuristic schemes such as local search with high-level diversification strategies. A large part of metaheuristics are built on top of a basic \emph{neighborhood search} (NS) scheme, where NS is a local search procedure that starts from an initial solution $\xx$ and iteratively search for improving solutions by exploring a series of neighborhoods.



We consider an instance $p \in \mathcal{P}$ of a CO problem, where $\mathcal{P}$ is the set of problem instances. $\mathcal{X}$ is the solution space, i.e., the set of feasible solutions. The neighborhood is defined as follows. Let $\xx$ be a solution of an instance $p$ such that $\xx\in \mathcal{X}$. The \emph{neighborhood} $N(\xx)$ of solution $\xx$ is a subset of solutions defined from $\xx$ in the solution space $\mathcal{X}$, i.e., $N(\xx) \subseteq \mathcal{X}$.


In general, the structure of the neighborhoods and how these neighborhoods are explored, are designed according to the characteristics of the problem at hand. A neighborhood structure is typically determined by a \emph{transformation} operator that applies a move to the current solution in the solution space.

A transformation operator is typically represented by $\Delta : \mathcal{X} \to  2^{\mathcal{X}}$, which is a function (or a set of functions) that maps a solution $\xx$ to a set of solutions. If the operator $\Delta$ is defined with a set $M$ of parameters, the operator can be defined  by $\Delta : \mathcal{X} \times M \to  2^{\mathcal{X}}$. For instance, a neighborhood $N(\xx) \mapsto \Delta (\xx, \mm)$ can be constructed by applying the transformation operator with $\mm \in M$ to the current solution $\xx$.


At each NS iteration, the neighborhood $N(\xx)$ is generated by a $ \Delta(\xx)$ or  
$\Delta (\xx, \mm)$ 
function.
A basic template for NS-based metaheuristic is shown in Algorithm \ref{alg:the_alg_meta}.

\begin{algorithm}[t]
\caption{NS-based metaheuristic}
\label{alg:the_alg_meta}
\SetAlgoLined
\KwIn{an initial solution $\xx$;}

$\xx^*\gets \xx$;

\Repeat {termination condition is reached} 
    {
    $N(\xx) \gets \Delta(\xx)$ \;\;(or $N(\xx) \gets \Delta(\xx,\mm)$ with $\mm \in M$);
    
    $\xx' \gets \argmin_{\xx'' \in N(\xx)} {f(\xx'')}$;
    
    $\xx \gets \xx';$
    
    \If {$f(\xx) < f(\xx^*)$} {$\xx^* \gets \xx;$}
    
    }

\Return $\xx^*$
\end{algorithm}

\subsection{Representation Learning for CO}
\label{sec:background_ml}
In ML, there is a vast library of models for representing CO problems depending on the format of input data of the CO task. For example, the model could be a linear function or some non-linear Artificial Neural Network (ANN) with a set of parameters to be optimized. In Deep Learning (DL), there are many types of neural networks and architectures available for modeling CO problems. For instance, a Multilayer Perceptron (MLP) is the simplest architecture of feedforward neural networks and can be used to model problems with fixed-size input data; Graph Neural Networks (GNNs) are developed to process  data that are naturally representable with graphs; Recurrent Neural Networks (RNNs) can be applied to process CO problems with sequential data, etc.

In particular, due to the ubiquity of graph data, problems over graphs arise in numerous application domains. Moreover, given the fact that the vast majority of CO problems have a discrete nature, many of them are naturally described in graphs or can be modeled into a graph structure. For instance, a network optimization problem can be naturally modeled by a graph and a generic MIP instance can also be represented into a bipartite graph \cite{gasse2019exact}. These graphs have inherent structural commonalities and patterns, and there is a potential to exploit valuable graph features to learn patterns from data. Therefore, \emph{graph representation learning} with the application of various GNNs \cite{gori2005new,scarselli2008graph,hamilton2017representation,cappart2021combinatorial} has recently emerged as a popular approach for studying CO with a machine learning perspective. Without loss of generality, GNNs learn a graph embedding based on the input features and the graph structure. In a nutshell, higher-level representations of a node or an edge are obtained by kernel convolutions which leverage its local structure.

In the literature, there has been an effort to learn algorithms for solving specific CO problems and many of them apply GNNs as the representation model. On the one hand, the first attempted paradigm was ``end-to-end'' learning for generating a heuristic solution by a ML model \cite{dai2017learning,nazari2018reinforcement, zhang2020learning,bello2016neural,gao2020learn,liu2021learning}.
However, these methods are typically limited to specific CO problems in which a heuristic solution can be easily constructed, and scaling to large-size instances is an issue. On the other hand, since a wide range of constrained CO problems can be formulated into a MIP model, there has also been increasing interest in learning decision rules to improve MIP algorithms \cite{gasse2019exact,he2014learning,Khalil_LeBodic_Song_Nemhauser_Dilkina_2016,khalil2017learning,balcan2018learning,Liu_Fischetti_Lodi_2022}. While it is shown that this direction has a great potential to improve the state-of-the-art of MIP algorithms, convincing generalization performances, and transfer learning across instances have not been fully tackled yet.

\section{Methodology}
\label{sec:methodology}

In this section, we present our framework for learning to generate high-quality solution neighbors for MH search. In Section \ref{sec:solution_structure} and \ref{sec:structural}, we will first introduce an efficient neighborhood structure by applying a transformation to a selected subset of variables from the current solution. Then, in Section \ref{sec:learning} the problem becomes to determine how to select a subset of variables that leads to a promising neighborhood, i.e., one containing a descent direction of the objective in the solution space. By conducting a classification task, our method learns promising variable selection policies from the structural characteristics of the problem and high-quality solutions. 


\subsection{Transformation Operator and Neighbor Generation}
\label{sec:solution_structure}

The definition of the neighborhood plays a critical role in the  NS-based metaheuristics. It specifies how the metaheuristic search moves in the solution space. As introduced in Section \ref{sec: background_mh}, the structure of a neighborhood is typically determined by a transformation operator, since the latter specifies how the current solution is perturbed and transformed to solution neighbors. 

For instance, in a classical Traveling Salesman Problem (TSP) where a set of $n$ cities and the distances between each pair of cities are given, the task is to find the shortest tour that visits each city exactly once. Given an arbitrary tour as the initial solution, a simple transformation operator can be defined by firstly removing $k$ edges from the current tour and then adding $k$ other edges to construct a new tour. This is known as the ``$k$-OPT'' neighborhood.  For a generic CO problem at hand, there are many ways of defining a transformation operator. In general, two main aspects must be taken into consideration: \emph{exploitation} (or \emph{intensification}) and \emph{exploration} (or \emph{diversification}) of the solution space.

On the one hand, the smaller the perturbation induced by the operator, the closer the constructed neighborhood to the current solution. The metaheuristic search will thoroughly exploit the local solution regions around the current solution. On the other hand, with more perturbation or more randomness induced by the operator, the neighbors can be defined far from the current solution. The heuristic search will have a larger chance to explore more solution regions that have been less visited before. On the extreme case, when an operator completely perturbs the solution or the solution is allowed to be fully changed, the neighborhood could be expanded to the entire solution space, and the complexity of exploring this neighborhood will be very high, as high as solving the original problem. Instead, the heuristic search will become totally random when the operator perturbs a part of the solution randomly. Moreover, the design of the transformation operator also strongly depends on the type of problem and its representation. The effectiveness of an operator might not be the same on different types of problems.

In order to guide the metaheuristic search to promising regions of the solution space, we will present a general ML framework for neighbor generation. More precisely, we will design transformation operators with parameterizations that are able to construct neighborhoods containing high-quality solutions.


\subsection{Variable Selection for Neighbor Generation}
\label{sec:structural}
In this section, we will define a neighborhood structure in NS-based metaheuristics. As mentioned before, we aim at generating promising neighbors by efficient transformation operators. 

In NS, improving solutions are typically found from neighbors defined by perturbing only a part of the solution. In practice, in order to control the size of the neighborhood, the ratio of perturbation is typically set to a relatively small value.
Although this ratio can be increased by diversification strategies when no improving solution is found from the last search, it is still very common that the best solutions found by two consecutive local search iterations share a large part of variable values. Nevertheless, this ``partial evolution'' of the local optima induces a class of transformation operators for constructing structural neighborhoods, where the transformation is defined on a subset of variables.

\begin{definition}
\label{def:operator_2}
A``subset'' transformation operator is a function $\Delta : \mathcal{X} \times \mathcal{Y} \to 2 ^{\mathcal{X}}$, where $\mathcal{X}$ denotes the solution space, $\mathcal{V}$ is the index set of variables, and $\mathcal{Y} = \{ 0, \;1\}^ {|\mathcal{V}|}$ defines a binary decision space for selecting a subset of $\mathcal{V}$. Consequently, for $\yy \in \mathcal{Y}$, the neighborhood $N(\xx)$ of a solution $\xx \in \mathcal{X}$ is defined as $\Delta (\xx, \yy)$, and only variables $x_i$ with $y_i=1$ can have their value changed by transforming $\xx$ into a neighbor in $N(\xx)$.
\end{definition}



When the components of $\yy$ are all set to $1$, all the variables in $\mathcal{V}$ will be allowed to change. As discussed before, generating neighbors over the entire variable set $\mathcal{V}$ might result in a large local problem. On the other hand, the ``partial evolution'' of improving solutions in local search also indicates that it is possible to define efficient transformation operators only on a subset of variables. Moreover, common characteristics are often present in good solutions in many applications. For instance, in a classical TSP problem, two cities that are far from each other are typically disconnected in good solutions, whereas cities close to each other have a higher probability of being connected. Hence, ML techniques can be employed to learn structural information from high-quality solutions and select a subset of variables that has a high probability of defining a neighborhood that contains a descent point of the objective in the solution space.


\subsection{Learning a Variable Selection Policy for Structural Neighbor Generation}
\label{sec:learning}



Given an instance of a generic CO problem and a solution $\xx \in \mathcal{X}$, we denote its current state as $\ss \in \mathcal{S}$, where $\mathcal{S}$ is the state space of the problem and typically consists of a set of selected features including $\xx$.
In order to define a structural neighborhood using the subset transformation operator (defined by Definition \ref{def:operator_2}), an instantiation of $\mathcal{Y}$ is required to select a subset of variables. Hence, a variable selection policy $\pi$ for selecting $\yy \in \mathcal{Y}$ can be defined by
\begin{align*}
    \pi : \mathcal{S} & \longrightarrow \mathcal{Y}\\
         \ss &\longmapsto \yy.
\end{align*}
Now, the question is
\begin{quote}
How to design a variable selection policy by which the subset transformation operator defines neighborhoods containing high-quality solutions?
\end{quote}

As mentioned before, the high-quality solutions found during the metaheuristic search often share common characteristics and patterns, and more importantly, many improved solutions in the neighborhood often share a part of the variables with the same values. Hence, we propose a general learning-based framework to extract these characteristics from data, and exploit the learned knowledge to guide the metaheuristic search towards compact and promising solution regions. Specifically, the framework learns a variable selection policy and the high-quality neighborhood structure will then be defined by selecting promising variables, the values of which are allowed to be changed from the current solution. 

\begin{figure}[H]
          \centering
          \includegraphics[width=0.8\linewidth]{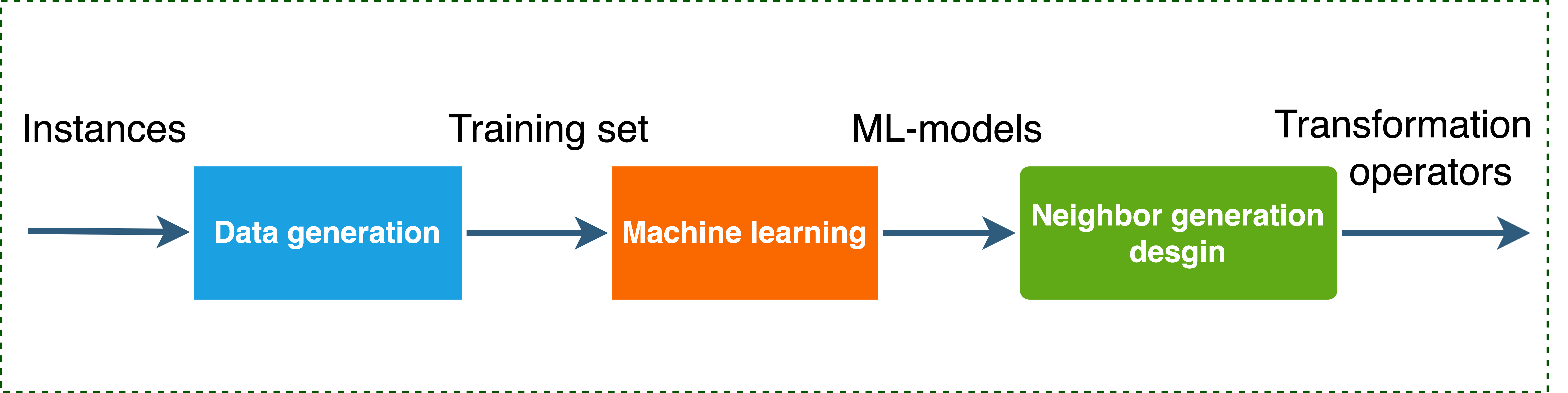}
          \caption{A learning-based framework for neighbor generation.}
          \label{fig:plot_tlb}
\end{figure}

The pipeline of the framework consists of three components: \emph{data generation}, \emph{machine learning} and \emph{neighbor generation design}. The framework is depicted in Figure \ref{fig:plot_tlb}.
 
In the following, we will instantiate the learning problem for variable selection as a classification task and train a variable selection policy through supervised learning. Within our method, the solution neighbors in metaheuristic search are generated based on our trained model. 



\subsubsection{Generation of a Training Data Set}
\label{sec:learn_datagene}

We aim to learn a policy $\pi$ that maps the state $\ss$ of a problem instance, to a high-quality $\yy$ (referred to as \emph{label} in the following), the binary classification decisions for selecting the best subset of variables that leads to a successful change of the current solution. For this purpose, one typically generate a training data set by applying 
an algorithm (referred to as \emph{expert}) which explores the neighborhood of some states $\ss_1,\ldots,\ss_N$, allowing to modify the values of all the variables, and then detects which are the variables whose value modification results in an improvement. Hence a label $\yy_i$ is associated with each $\ss_i$ ($i=1,
\ldots,N$) and we thus have a training data set $\mathcal{D}_{\text{train}} = \{(\ss_i,\yy_i))\}_{i=1}^N$.

\subsubsection{Training the Variable Selection Policy}
\label{sec:learn_train}
The policy we want to determine depends on parameters $\theta \in \Theta$ of a model (see below), where $\Theta$ is the set of candidate parameters. Each policy is thus denoted $\pi_{\theta}$ and, given a training data set $\mathcal{D}_{\text{train}} = \{(\ss_i,\yy_i))\}_{i=1}^N$, we aim to determine parameters $\theta^*$ such that $\pi_{\theta^*}(\ss_i)$ is as close as possible to $\yy_i$ ($i=1,\ldots,N)$. For this purpose, we solve the following classification problem:
\begin{equation}
\label{eq:training}
    \theta^* = \argmin_{\theta \in \Theta} \sum_{i=1}^N \mathcal{L} \left(\pi_{\theta}(\ss_{i}), \yy_{i} \right),
\end{equation}
where $\mathcal{L} \left(\pi_{\theta}(\ss_{i}), \yy_{i} \right)$ denotes the loss function and will be defined according to the CO application (see more details in Sections \ref{sec:app1_learn_drop}, \ref{sec:app1_learn_add} and \ref{sec:app2_learn}). In summary, given a state $\ss\in \mathcal{S}$, we use $\yy=\pi_{\theta^*}(\ss)$ for selecting which variables are allowed to be changed when moving to neighbors.


\paragraph{Modeling}

 Since CO problems generally have a discrete nature and many of them can be represented by a graph structure, we will apply GNN architectures \cite{gori2005new,scarselli2008graph, hamilton2017representation} for parametrizing the variable selection policy.
More precisely, if the state of a CO instance is modeled by a graph, GNNs can be implemented as the representation model for $\pi_{\theta}$. GNNs exhibit some appealing properties for processing data in graph format. First of all, GNNs are size-and-order invariant to input data, i.e., they can process graphs of arbitrary size and topology, and the graph model is invariant to the ordering of the input elements, which brings a critical advantage compared to other neural networks. Moreover, GNNs can exploit the sparsity of the graph by localized propagation of information, making them an ideal class of models for embedding sparse CO problems \cite{gasse2019exact}.



The basic architecture of GNNs consists of 3 modules: the input module, the convolution module, and the output module. In the input module, the state $\ss$ of the problem is fed into the GNN model. The input module embeds the features of $\ss$. The convolution module propagates the features of the graph components by \emph{graph convolution layers} \cite{gasse2019exact, Liu_Fischetti_Lodi_2022} and further embeds the features. 

In particular, the architecture of the graph convolution layers used in this paper applies the message passing operator, defined as
\begin{align}
    \mathbf{v}_i^{(h)} = f_{\psi}^{(h)} \left( \mathbf{v}_i^{(h-1)}, \sum_{j \in \mathcal{N}(i)} \, g_{\phi}^{(h)}\left(\mathbf{v}_i^{(h-1)}, \mathbf{v}_j^{(h-1)},\mathbf{e}_{j,i}\right) \right),
\end{align}
where $\mathbf{v}^{(h-1)}_i \in \mathbb{R}^d$ denotes the feature vector of node $i$ in layer $(h-1)$, $\mathbf{e}_{j,i} \in \mathbb{R}^m$ denotes the feature vector of edge $(j, i)$ from node $j$ to node $i$, $\mathcal{N}(i)$ denotes the set of nodes adjacent to $i$ in the graph, and $f_{\psi}^{(h)}$ and $ g_{\phi}^{(h)}$ denote the embedding functions paramaterized by $\psi$ and $\phi$ in layer $h$ and are typically represented by neural networks (e.g., MLP architectures).
The output module maps the embedding of the state into a distribution in the binary decision space for each variable and the final output is the class prediction of the variable (typically 1 for being selected into the subset, 0 for not being selected).

\subsubsection{Neighbor Generation Design} 
\label{sec:learn_neigh}
After training the classification model for variable selection, the next step is to apply the pretrained model for neighbor generation. It is important to note that, the trained classification model itself is probabilistic, and it maps the current solution state of an instance into a probability distribution in the binary decision space for each variable. One still needs to select a strategy to make binary decisions on the variables. 

The most straightforward way is to apply a \emph{greedy} strategy. The decision is made by always picking the class with a higher probability and the resulting strategy is \emph{deterministic}. Another way is to sample decisions from the distribution, and hence the strategy is \emph{probabilistic}. In general, the greedy strategy only selects the deterministic subset of variables and exploits the learned knowledge from the training data set. Instead, the probabilistic strategy selects from all the possible neighborhoods with a probability preference defined by the classification model, thus results in a better exploration of the less visited solution regions. There is no guarantee that one strategy is always better than the other. In practice, they can be combined. For each task, one should search for a good trade-off between exploitation of the local solution regions and exploration of the global solution space. 

Finally, the template for a NS-based metaheuristic guided by the variable selection policy $\pi$ is given in Algorithm \ref{alg:the_alg_meta_ml}. The neighborhood structure $N(\xx)$ is defined by the subset transformation operator that selects a subset of variables according to policy $\pi_{\theta^*}$, and modifies some of their values to lead the heuristic search to promising solution regions.  

\begin{algorithm}[t]
\caption{NS-based metaheuristic guided by the variable selection policy}
\label{alg:the_alg_meta_ml}
\SetAlgoLined
\KwIn{an initial solution $\xx$ and  a variable selection policy ${\pi_{\theta^*}}$;}
$\xx^* \gets \xx$;

\Repeat {termination condition is reached} 
    {
    $\yy \gets \pi_{\theta^*}(\ss)$ where $\ss$ is the state associated with $\xx$;
    
    $N({\xx}) \gets \Delta({\xx},\yy)$;
    
    ${\xx}' \gets \argmin_{{\xx}'' \in N({\xx})} {f({\xx}'')}$;
    
    ${\xx} \gets {\xx}'$;
    
    \If {$f({\xx}) < f({\xx^*})$} { 
         ${\xx^*} \gets {\xx}$;}

    }

\Return ${\xx^*}$
\end{algorithm}

\section{Application 1: Tabu Search in Wireless Network Optimization}
\label{sec:app1}

In this section, we apply our framework to the first case study, a Wireless Network Optimization (WNO) problem. We will demonstrate the effectiveness of our learning-based framework by generating solution neighbors in a Tabu Search scheme. 


\subsection{The Tactical WNO Problem}

When telecommunications are necessary but standard networks are unavailable, as in the case of disaster relief operations, small temporary wireless networks, called tactical networks, are set up. These typically connect between 10 and 50 nodes (the key locations that must communicate) in a single network. The design of such networks can be optimized such that the network's weakest link is maximized. This, in turn, guarantees that all nodes can receive important information in a timely manner.

Tactical wireless network design is a complex non-linear combinatorial optimization problem that includes three sub-problems: the design of the topology ($P_0$), the configuration of the network ($P_1$), and the configuration of the antennas ($P_2$). These sub-problems are nested such that $P_1$ is defined given a topology, and $P_2$ is defined given a topology and a network configuration, namely
$$
\underbrace{\max_{\substack{\textrm{topology }t\in T\\~\\~}}}_{P_0} \ \ \underbrace{\max_{\substack{\textrm{network}\\\textrm{configuration}\\~}}}_{P_1} \ \ \underbrace{\max_{\substack{\textrm{antenna}\\\textrm{configurations}\\~}} o.}_{P_2}
$$
where $T$ is the space of valid topologies and $o$ is the objective function of the problem. The problem has been proposed by an industrial partner and the modeling of the full tactical wireless network design optimization problem can be found in \cite{perreault2022}.


Given a set of nodes $V$, the topology $t \in T$ can be any undirected tree $(V,E)$ that describes how information travels in the network between every pair of nodes. Each edge in the topology represents a direct connection between two antennas in the network. A network configuration selects a root node (also known as a master hub or gateway) and assigns waveforms and channels/frequencies to the edges. The waveforms are the communication protocols that depend on the local structure of the edges and the channels and frequencies characterize their radio signals. Together with the antenna configurations, they determine which edges interfere with each other. An antenna configuration requires to define an angular alignment with respect to the azimuth as well as a set of activated beams, in the case of multi-beam antennas. Given all these properties, the radio signals can be physically modeled by also taking into account the path losses and fade margins of the terrain between every pair of nodes. The data transmission speed (called direct throughput) $T\!P_{uv}$ for all the edges $[u,v] \in E$ can then be computed and, depending on the traffic scenario $X$ and its distribution of congestion $n_{uv}^X$ in the edges, their effective throughputs $T\!P_{uv}/n_{uv}^X$ can also be computed. There are three traffic scenarios: scenario $A$ is a single communication between any two nodes, scenario $B$ is all nodes communicating with the root node, and scenario $C$ is all nodes communicating with each other. For all traffic scenarios, the congestion of any directed edge $(u,v)$ can be computed as a function of the number of descendants $\textrm{desc}_v$ defined by the root node selected in the network configuration.

To evaluate a single topology $t \in T$, problem $P_1$, which is itself a complex combinatorial optimization problem, must be solved, and solving it requires solving $P_2$ many times as well. In turn, problem $P_2$ can be efficiently solved by a simple geometrically-based heuristic. Solving $P_1$ to optimality every time that a topology needs to be evaluated is too costly, especially in a NS context. It can be approximated in such a way that it can be solved efficiently by exhaustive enumeration. We denote the resulting approximated objective function by $\overline{f}(t)$. However, this approximation does not necessary provide a feasible network because it does not define a proper frequency assignment. Alternatively, $P_1$ can also be tackled directly by considering greedy frequency assignments which do provide feasible networks, although this takes considerably more time. We denote this final objective function by $f(t)$.

Given such procedures for solving $P_1$ and $P_2$, $P_0$, given by
$$
\max_{t\in T} f(t),
$$
can be solved with a NS-based MH in the space of topologies, where local directions of descent are computed with $\overline{f}(t)$.

\subsection{Topology Tabu Search}

\subsubsection{Neighborhoods}

The neighborhood structure dictates how the local search can move in the space of solutions. 
Our neighborhood $N$ is the edge-swap neighborhood. To construct the neighbor topologies, for each tree edge $e \in E$ in the current topology, we remove it, which disconnects the topology into two connected components, and we consider every possible way of reconnecting these two connected components with another edge. Hence, $\Delta(t)$ is the set of all topologies that can be obtained from $t$ by  an edge-swap move.

\begin{figure}[ht]
\centering
\parbox{1.5in}{
  \centering
  \includegraphics[width=1.5in]{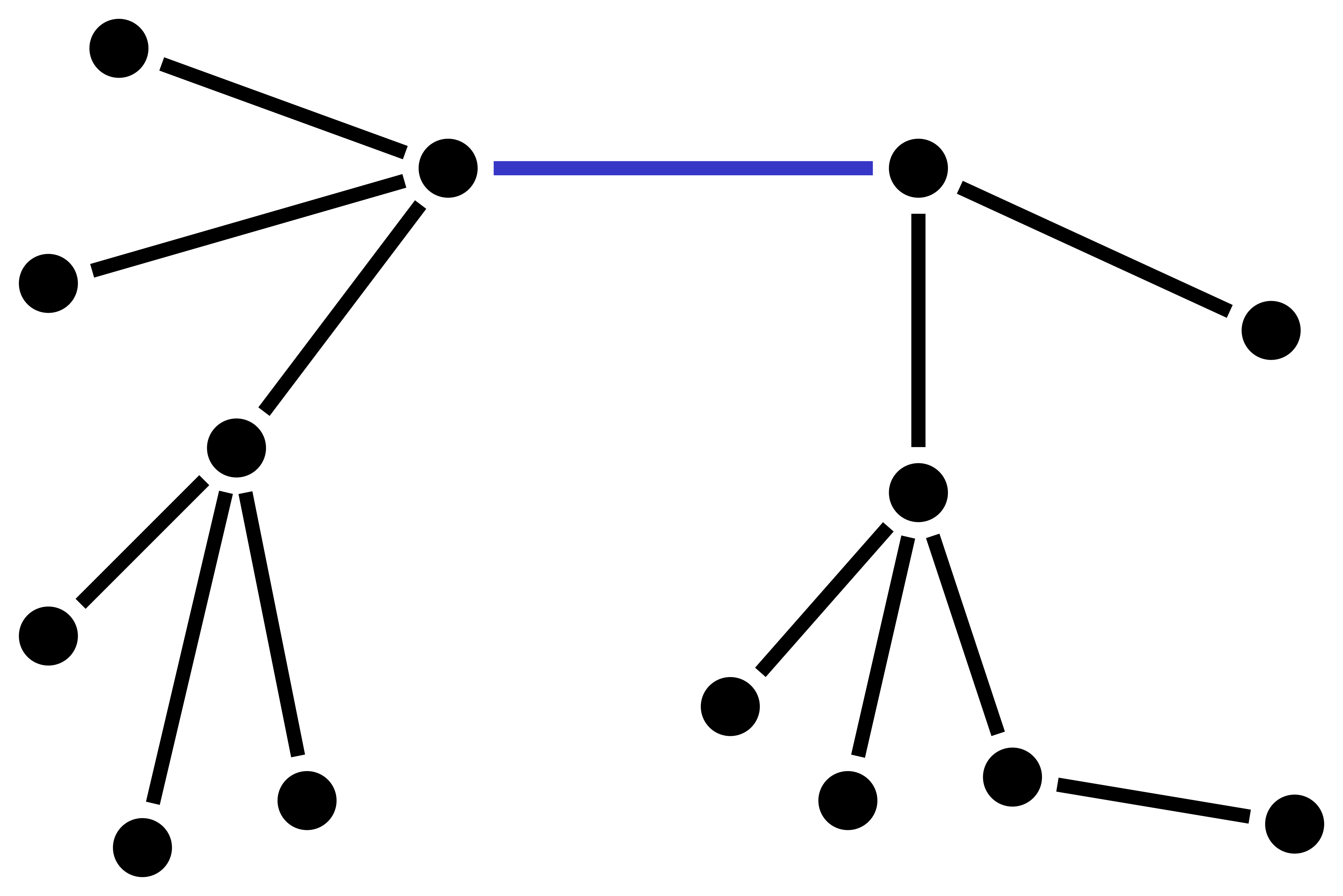}
}
\ \ \ \scalebox{1.05}{\contour{black}{$\leftrightarrow$}} \ \ \
\parbox{1.5in}{
  \centering
  \includegraphics[width=1.5in]{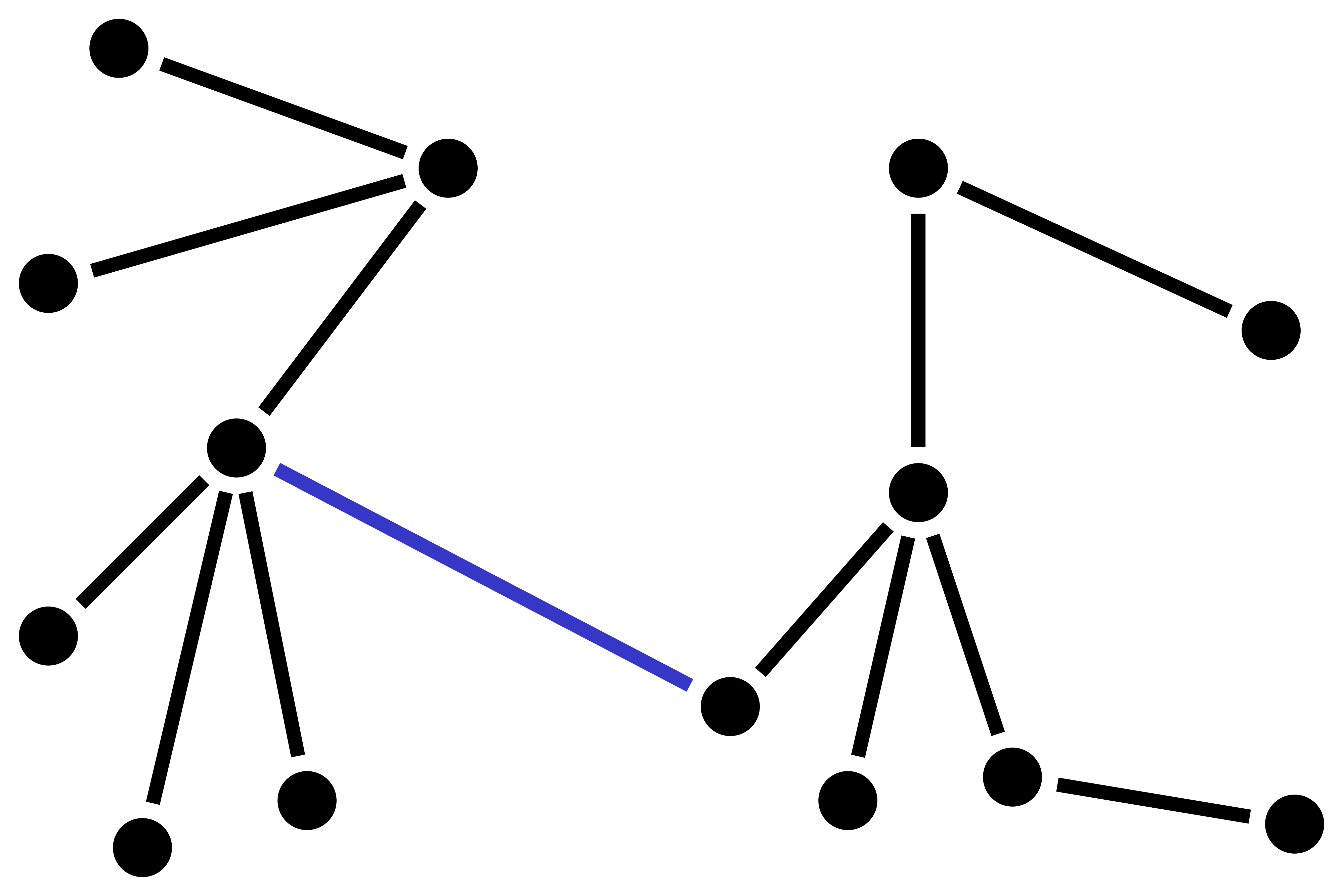}
}
\caption{Neighbor Topologies.} \label{fig:topology}
\end{figure}

An example of neighbor topologies is depicted in Figure \ref{fig:topology}, where the edges in blue were ``swapped''.

\subsubsection{Tabu Search}

Tabu Search (TS) is a MH that explores the solution space by temporarily not allowing to move in the direction from which it came. The only case in which this rule does not apply is when moving in such a direction improves the incumbent solution. This is done by using tabu lists that record the opposites of the last few moves, which then become tabu for the next few moves.

For our edge-swap neighborhood, two tabu lists are necessary:  a list $L_{drop}$ for forbidding edges to be dropped and a list $L_{add}$ for forbidding edges to be added. Let $t'\in N(t)$ be a neighbor of the current topology $t$, obtained by replacing an edge $e$ with an edge $e'$ : we consider $t'$ as tabu if $e\in L_{drop}$ or (non exclusive) $e'\in L_{add}$. When a move is actually made from the current topology, the newly dropped edge is added to ${L_{add}}$ (it cannot be added back for the next few moves) and the newly added edge is added to ${L_{drop}}$ (it can not be dropped for the next few moves).

The lengths of the tabu lists determine for how many moves the dropped/added edges have tabu status. As observed in \cite{perreault2022},  reasonable lengths for ${L_{drop}}$ and ${L_{add}}$  are
$
    \left\lfloor (\sqrt{|V|-1})/2\right\rceil
$
and 
$
    \left\lfloor\sqrt{|V| \, (|V| - 1)/2}\right\rceil
$, respectively,  
where $\left\lfloor x \right\rceil$ means $x$ rounded to the nearest integer. The pseudocode for the resulting \emph{topology Tabu Search} metaheuristic is given in Algorithm \ref{alg:P0tabu}.

\begin{algorithm}[t]
\caption{$P_0$ Topology Tabu Search}
\label{alg:P0tabu}
\SetAlgoLined
\KwIn{an initial topology $t$;}
$t^* \leftarrow t$; $\overline{t}^* \leftarrow t$; $L_{drop} \leftarrow \emptyset$; $L_{add}\leftarrow \emptyset$;

\Repeat {termination condition is reached}
{
    $N(t) \gets \Delta(t)$;
    
    $t'\gets \argmax_{t'' \in N(t)}\overline{f}(t'')$\;
 
   \hspace{1cm}$\textrm{s.t. }  t'' \textrm{ is not tabu or }  \overline{f}(t'') > \overline{f}(\overline{t}^*)$;
    
    \vspace{5pt}\If {$\overline{f}(t') > \overline{f}(t^*)$} { 
        $\overline{t}^* \leftarrow t'$;
    }
    
    Make move $t \leftarrow t'$ and update tabu lists $L_{drop}$ and $L_{add}$;
    
    \If {$f(t) > f(t^*)$} { 
        $t^* \leftarrow t$;
    }
}
\Return $t^*$
\end{algorithm}

\subsection{Learning to Generate Edge-Swap Neighbors for TS}
\label{sec:app1_learn}

As a NS-based MH, TS can be roughly described as a NS scheme plus a ``tabu'' strategy for preventing cycling by keeping a short-term memory of visited solutions stored in the tabu list. As discussed above, the topology design for the wireless network optimization problem can be solved by applying a topology TS algorithm (Algorithm \ref{alg:P0tabu}). In this TS scheme, a neighborhood $N(t)$ of the current solution $t$ is constructed by applying an edge-swap move operator. Specifically, the move operator consists of two steps: an edge will be dropped from the current topology (by enforcing the value of the dropped edge variable from $1$ to $0$) and another edge will be added to complete a new topology (by enforcing the value of the added edge variable from $0$ to $1$). As a result, the edge-swap neighborhood consists of all the edge-swap moves. 

In Algorithm \ref{alg:P0tabu}, the neighborhood search at each TS iteration explores the entire edge-swap neighborhood by enumerating all  possible moves. Although the algorithm is guaranteed to find a locally optimal solution by \emph{exploiting} the entire edge-swap neighborhood, it is generally slow in terms of computing time since a complex subproblem ($P_1$) has to be solved to evaluate each ``edge-swap'' move, and another subproblem ($P_2$) has to be solved to evaluate each solution in $P_1$. 

\begin{figure}[ht]
          \centering
          \includegraphics[width=0.8\linewidth]{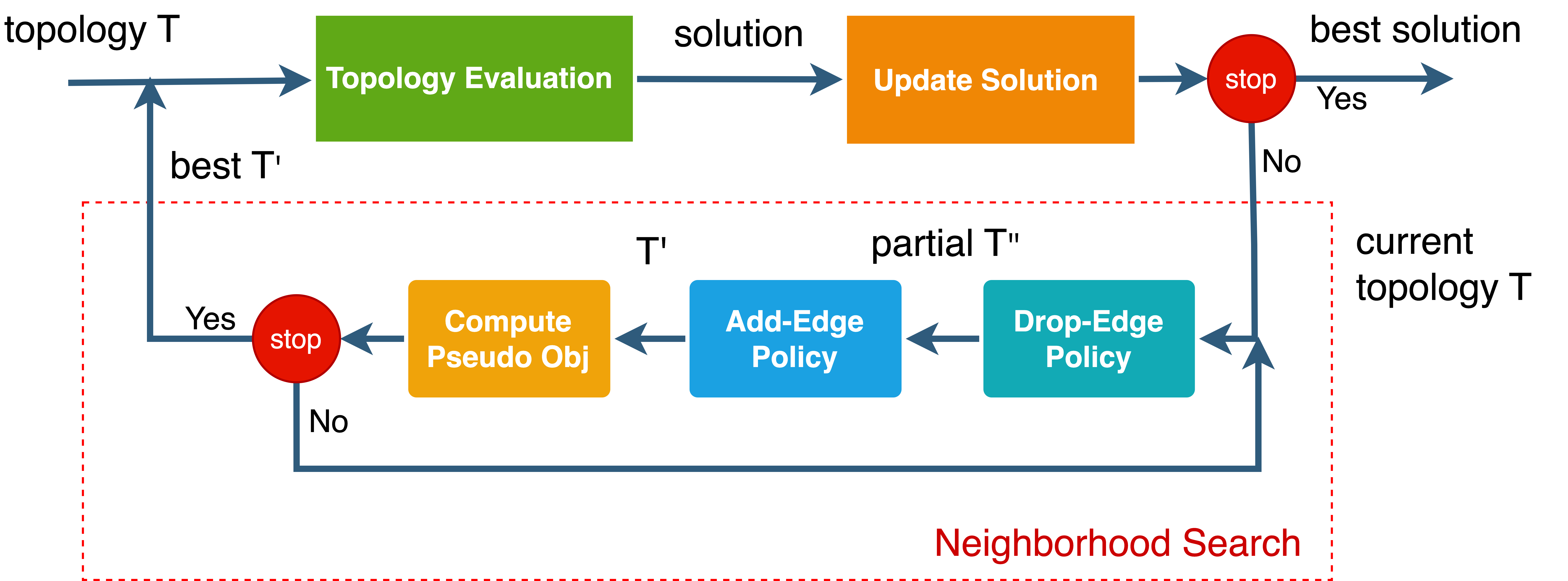}
          \caption{A learning-based framework for generating ``edge-swap'' neighbors.}
          \label{fig:plot_tlb_2}
\end{figure}

A potential improvement of the enumeration strategy is to apply a random strategy to sample from droppable edge variables and addable edge variables, thus evaluating only a subset of possible moves. The random strategy generally \emph{explores} a larger solution space than the enumeration strategy because, if the overall amount of time for the algorithm is fixed, it is able to do more iterations. However, random sampling might not be efficient enough for guiding the search towards high-quality neighborhoods. To achieve a better trade-off between exploitation and exploration, we propose to exploit structural characteristics of improving moves, and learn good variable selection policies for selecting ``drop'' edge variables and ``add'' edge variables to generate size-reduced, but high-quality edge-swap neighborhoods.
The scheme for variable selection consists of two components: a classifier for selecting the ``drop'' edge variables and another classifier for selecting the ``add'' edge variables. The scheme is depicted in Figure \ref{fig:plot_tlb_2} and a detailed pseudocode will be given in Algorithm \ref{alg:P0tabuGNN}.
 
\subsubsection{Learning to Drop Edges}
\label{sec:app1_learn_drop}

Let $\mathcal{S}$ denote the state space of the topologies of a WNO instance, and let $\mathcal{Y}^d = \{0,1\}^{|\mathcal{I}|}$ denote the set of possible binary decisions for the droppable edge variables, where $\mathcal{I}$ is the index set of droppable edges in the current topology. We aim to learn a ``drop-edge'' policy that maps each state $\ss\in \mathcal{S}$ to a label $\pi^d_{\theta^*}(\ss)=\yy^d\in \mathcal{Y}^d$. 
The generation of a training dataset is explained in Section \ref{sec:app1_data}. Roughly speaking, we consider a set $J$ of WNO instances, and apply $T_j$ iterations of Algorithm \ref{alg:P0tabu} to each instance $j\in J$. At each iteration of each instance, we determine the edges that can be dropped and replaced by another edge to produce a better topology. This gives a label for every visited state. Let 
$\{\ss_{i}^{(j)}\}_{i=1}^{T_j}$ be the states of the $j^{th}$ instance, and let    $\{\yy_{i}^{d(j)}\}_{i=1}^{T_j}$ be the corresponding labels. As explained in Section \ref{sec:learn_train}, the policy $\pi^d_{\theta^*}$ is obtained by solving the following problem: 
\begin{equation}\label{eq:maxlikelihood}
    \theta^* = \argmin_{\theta \in \Theta} \sum_{j=1}^{|J|} \sum_{i=1}^{T_j} \mathcal{L} \left(\pi^d_{\theta}(\ss_{i}^{(j)}), \yy_{i}^{d(j)} \right),
\end{equation}
where $\mathcal{L}$ is the loss function.

\paragraph{Feature Design}
Since the topologes naturally have a graph structure, we represent each state $\ss$ by a graph. More precisely, given a WNO instance with $n$ nodes, we consider $d$ features for each node and $e$ features for each edge. The matrix of node features is denoted by $\mathbf{V} \in \mathbb{R}^{n\times d}$ and the tensor of edge features is denoted by $\mathbf{E} \in \mathbb{R}^{n\times n\times e}$. The features used for learning which edges  have the best potential to lead to an improvement when dropped are listed in Table \ref{table:drop-features}. 

\begin{table}[ht]
    \centering
    \begin{tabular}{c l l}
    \multicolumn{1}{c}{Tensor} & \multicolumn{1}{l}{Feature} & \multicolumn{1}{l}{Description} \\
    \toprule
    \multirow{1}{*}{$\mathbf{V}$} 
    & coordinates & $X$ and $Y$ coordinates in km \\
     \cmidrule{2-3}
    & desc\_v & normalized $desc_v = ( (|desc_v|-1)/(|V|-1))$ \\
    \midrule
    \multirow{1}{*}{$\mathbf{E}$}
    & pass\_loss & path loss in dB \\
    \cmidrule{2-3}
    & fade\_margin & fade margin in dB \\
    \cmidrule{2-3}
    & waveform & binary (point-to-point or point-to-multipoint) \\
    \cmidrule{2-3}
    & channel & frequency channel \\
    \cmidrule{2-3}
    & n\_beams & number of beams used by the antenna of the predecessor\\
    \cmidrule{2-3}
    & tp\_uv\_a & throughput scenario A\\
    \cmidrule{2-3}
    & tp\_uv\_b & throughput scenario B\\
    \cmidrule{2-3}
    & tp\_uv\_c & throughput scenario C\\
    \cmidrule{2-3}
    & tp\_uv\_m & throughput by mixed scenarios of A, B, C\\
    \cmidrule{2-3}
    & edge indices & edge indices in current topology \\
    \bottomrule
    \end{tabular}
    \caption{Description of the features in the ``droppable'' graph $(\mathbf{N}, \mathbf{E})$. More detailed definitions of the listed features can be found in \cite{perreault2022}.}
    \vspace{10pt}
    \label{table:drop-features}
\end{table}

\paragraph{GNN Model}
Because of the intrinsic graph structure of the topologies, it is natural to apply GNN to model the drop-edge classifier.
The chosen GNN architecture has 3 modules: the input module, the convolution module and the output module. The output module embeds the hidden features extracted from the convolution module and maps the embedding of each edge in the current topology into a two-neuron output. 

\paragraph{Loss Function}
From the data generation process, we observed that class distribution is unbalanced and only less than $50\%$ of droppable edges lead to an improving move. Since the objective of the learning is to select as many ``improving'' edges to be dropped as possible, it is reasonable to put more effort in improving the predictions on the minority class. Therefore, we apply a \emph{weighted} Cross Entropy (WCE) loss to train the model. Specifically, we add a penalty factor $\lambda$ for the ``improving" class and $1 - \lambda$ for the ``non-improving" class. Formally, the WCE loss is defined as

\begin{equation}
     \mathcal{L}_{wce} \left(\hat{\yy}, \yy \right) = 
     \frac{1}{|\mathcal{I}|} \sum_{i=1}^{|\mathcal{I}|} \left( -\lambda y_i \log \hat{y}_i - (1-\lambda)(1 -y_i) \log (1- \hat{y}_i)\right),
\end{equation}
where $\hat{\yy}$ denotes the prediction of probability given by the policy model, $\yy$ denotes a label, and $\mathcal{I}$ is the index set of droppable edges in the current topology. The penalty factor $\lambda \in [0.5, 1]$ imposes a larger loss for the ``improving'' class during training.

\subsubsection{Learning to Add Edges}
\label{sec:app1_learn_add}
After dropping an edge, the next step is to select an edge to be added from all addable edges. Let $\mathcal{S}$ denote the state space of all graphs obtained from a topology by dropping an edge, and let $\mathcal{Y}^a = \{0,1\}^{|\mathcal{I}|}$ be the set of possible binary decisions for addable edges, where $\mathcal{I}$ is the index set of addable edges. We aim to learn an ``add-edge'' policy $\pi^a_{\theta^*}$ that maps each state $\ss\in \mathcal{S}$ to a label $\yy^a\in 
\mathcal{Y}^a$. 

The generation of a a training dataset is similar to what was done for dropping edges. Algorithm \ref{alg:P0tabu} is applied to a set $J$ of WNO instances for a few iterations in order to determine which edges have the potential to improve a topology once one of their edges has been dropped. More details are given in Section \ref{sec:app1_data}.
Let 
$\{\ss_{i}^{(j)}\}_{i=1}^{T_j}$ be the states of the $j^{th}$ instance, and let    $\{\yy_{i}^{a(j)}\}_{i=1}^{T_j}$ be the corresponding labels. The policy $\pi^a_{\theta^*}$ is obtained by solving the following problem: 
\begin{equation}\label{eq:maxlikelihood_2}
    \theta^* = \argmin_{\theta \in \Theta} \sum_{j=1}^{|J|} \sum_{i=1}^{T_j} \mathcal{L} \left(\pi^a_{\theta}(\ss_{i}^{(j)}), \yy_{i}^{a(j)} \right),
\end{equation}.


\paragraph{Feature Design}
As was done for dropping edges, we model a state $\ss$ by a graph. Given a WNO instance with $n$ nodes, $d$ features for each node, and $e$ features for each edge, we denote by $\mathbf{V} \in \mathbb{R}^{n\times d}$ the matrix of node features and by $\mathbf{E} \in \mathbb{R}^{n\times n\times e}$ the tensor of edge features.
The features used for learning which edges have the best potential to lead to an improvement when added are listed in Table \ref{table:add-features}. 

\begin{table}[htbp!]
    \centering
    \begin{tabular}{c l l}
    \multicolumn{1}{c}{Tensor} & \multicolumn{1}{l}{Feature} & \multicolumn{1}{l}{Description} \\
    \toprule
    \multirow{1}{*}{$\mathbf{V}$} 
    & coordinates & $X$ and $Y$ coordinates in km \\
     \cmidrule{2-3}
    & desc\_v & normalized $|desc\_v| = ( (|desc\_v|-1)/(|V|-1))$ \\
    \midrule
    \multirow{1}{*}{$\mathbf{E}$}
    & edge\_type & binary indicator of addable edges \\
    \cmidrule{2-3}
    & pass\_loss & path loss in dB \\
    \cmidrule{2-3}
    & fade\_margin & fade margin in dB \\
    \cmidrule{2-3}
    & edge indices & edge indices in current topology \\
    \bottomrule
    \end{tabular}
    \caption{Description of the features in the ``addable'' graph $(\mathbf{N}, \mathbf{E})$.}
    \vspace{10pt}
    \label{table:add-features}
\end{table}

\paragraph{GNN Model}
As for the drop-edge classifier, we also apply GNN to model the add-edge classifier for selecting the edge to be added. The GNN architecture is the same as for the drop-edge classifier.

\paragraph{Loss Function}
The class distribution for addable edges is even more unbalanced than for the droppable case. Actually, only less than $10\%$ of addable edges lead to an improving move. Therefore, we also applied the WCE loss to train the model with an emphasis on the minority class.

\subsection{Numerical Experiments}
\label{sec: app1_experiment}

This section contains the experimental results for the WNO application. After presenting the data collection, in Sections \ref{sec:app1_setting} and \ref{sec:app1_metric}, we discuss the experimental setting and the evaluation metrics, respectively. Finally, the results are reported and discussed in Section \ref{sec:app1_results}.

\subsubsection{Data Collection}
\label{sec:app1_data}
\paragraph{Problem Instance Generation}
For the numerical experiments, the instances were generated in the following way. First, independent coordinates for the explicit nodes $v \in V$ are iteratively generated using
\begin{equation}
    (x_v,y_v) = \left( \sqrt{U_0} \cos (2 \pi U_1), \ \sqrt{U_2} \sin (2 \pi U_3)\right),
\end{equation}
where $U_0,U_1,U_2,U_3 \sim \mathcal{U}(0,1)$ are independent and identically distributed (IID). These coordinates are then scaled to match the average distance ratio of 10 km. This coordinate generation is repeated until the minimum distance is above 2 km and the maximum distance is below 150 km. Once this is achieved, for each possible edge $[u,v]$ with $u,v \in V$, a random uniform variable is sampled for the path loss and another is sampled for the fade margin, both according to empirical cumulative distribution functions derived from North-American datasets. 

We generated small instances of 10 nodes, as well as larger instances of 30 nodes. For these instances, the TS MH is initialized with a minimum spanning tree with a specific distance based on the path losses and fade margins of the terrain between every pair of nodes \cite{perreault2022}.

\paragraph{Training Data Generation}
To collect data for training the drop-edge and add-edge classifiers, we call Algorithm \ref{alg:P0tabu} to evaluate all possible moves in the edge-swap neighborhood, and compute labels according to the learning task. In particular, given a WNO instance, and an initialization of the topology, we execute Algorithm \ref{alg:P0tabu} as follows: for each TS iteration, we first call the TS algorithm to evaluate all possible edge-swap moves. For each droppable edge, we evaluate all the addable edges. If the best resulting edge-swap move within the neighborhood of the droppable edge leads to a better approximated pseudo-objective, then the droppable edge is labeled as ``improving''; otherwise, it is labeled as ``non-improving''. Given the dropped edge, the label of each addable edge is also decided by the quality of the resulting edge-swap move. If the move leads to a better approximated pseudo-objective, then the addable edge within the neigborhood of the corresponding dropped edge is labeled as ``improving''; otherwise, it is labeled as ``non-improving''.

\paragraph{TS Guided by GNNs}
The TS algorithm with GNN classifiers (TS-GNN) is obtained by using the learned policies $\pi^d_{\theta^*}$ and $\pi^a_{\theta^*}$ to generate edge-swap neighbors. Hence, the pseudocode of TS-GNN is  the same as Algorithm \ref{alg:P0tabu}, except that the neighborhood $N(t)$ of $t$ results from the subset transformation operator that performs edge-swaps according to $\pi^d_{\theta^*}$ and $\pi^a_{\theta^*}$. A detailed pseudocode is given in Algorithm \ref{alg:P0tabuGNN}.

\begin{algorithm}[t]
\caption{$P_0$ Topology Tabu Search with GNN classifiers}
\label{alg:P0tabuGNN}
\SetAlgoLined
\KwIn{an initial topology $t$, a variable selection policy $\pi_{\theta^*}^d$ for dropping edges,  and a variable selection policy $\pi_{\theta^*}^a$ for adding edges;
}
$t^* \leftarrow t$; $\overline{t}^* \leftarrow t$; $L_{drop} \leftarrow \emptyset$; $L_{add}\leftarrow \emptyset$;

\Repeat {termination condition is reached}
{
$N(t) \gets \emptyset$

$\yy^d \gets \pi_{\theta^*}^d(\ss)$ where $\ss$ is the state associated with $t$;

\For{every $e_i $ in $t$ such that $y_i^d = 1$}
    {
    
    set $\ss_i$ equal to the state obtained from $\ss$ by dropping $e_i$ from topology $t$;
    
    $\yy^a \gets \pi_{\theta^*}^a(\ss_i)$;

    \For{each reconnecting edge $e_j$ such that $y_j^a = 1$}
        {
        set $t_{i,j}$ equal to the topology obtained from $t$ by dropping edge $e_i$ and adding edge $e_j$;
        
        $N(t) \gets N(t) \cup \{t_{i,j}\}$;
    }
}
  
    $t'\gets \argmax_{t'' \in N(t)}\overline{f}(t'')$\;
 
   \hspace{1cm}$\textrm{s.t. }  t'' \textrm{ is not tabu or }  \overline{f}(t'') > \overline{f}(\overline{t}^*)$;
    
    \vspace{5pt}\If {$\overline{f}(t') > \overline{f}(t^*)$} { 
        $\overline{t}^* \leftarrow t'$;
    }
    
    Make move $t \leftarrow t'$ and update tabu lists $L_{drop}$ and $L_{add}$;
    
    \If {$f(t) > f(t^*)$} { 
        $t^* \leftarrow t$;
    }
}
\Return $t^*$
\end{algorithm}

\subsubsection{Experimental Setting}
\label{sec:app1_setting}

\paragraph{Training}
We train the drop-edge policy and add-edge policy on the collected training dataset separately. The training dataset is generated from 30 small instances of 10 nodes. The collected data is split into training ($70\%$), validation ($10\%$), and test ($20\%$) sets.

\paragraph{Evaluation}
We evaluated the compared algorithms listed in Section \ref{sec:app1_results} on two evaluation sets. The first evaluation set contains $50$ small instances of $10$ nodes. In addition, to evaluate the generalization performance of our approach on larger instances,  we also evaluate the model trained on small instances on a large evaluation set with $5$ instances of $30$ nodes.

\paragraph{Experimental Environment}
Our code is written in Python $3.9$ and we use Pytorch $1.60$ \cite{paszke2019pytorch}, Pytorch Geometric $1.7.0$ \cite{fey2019fast} for implementing and training the GNNs.

\subsubsection{Evaluation Metrics}
\label{sec:app1_metric}

We use the \emph{primal integral} \cite{berthold2013measuring} to measure the performance of the compared algorithms. The \emph{primal integral} was originally proposed to measure the performance of primal heuristics for solving mixed-integer programs. The metric takes into account both the quality of solutions and the computing time spent to find those solutions during the solving process. To define the primal integral, we first consider a \emph{primal gap function} $p(t)$ as a function of time, defined as
\begin{align*}
        p(t) = 
        \left\{
        \begin{array}{ll}
            1, & \ \ \ \text{if no incumbent until time $t$},\\
            \bar{\gamma} (\tilde{\xx}(t)), &\ \ \ \text{otherwise},
        \end{array}
        \right.
\end{align*}
where $\tilde{\xx}(t)$ is the incumbent solution at time $t$, and $\bar{\gamma} (\cdot) \in [0,1]$ is the \emph{scaled primal gap}
\begin{align*}
        \bar{\gamma} (\tilde{\xx}) = 
        \frac{| f(\tilde{\xx}_{\text{opt}}) - f(\tilde{\xx}) |}{| f(\tilde{\xx}_{\text{opt}}) - f(\tilde{\xx}_{\text{init}}) |},
\end{align*}
where $f(\tilde{\xx})$ denotes the objective value given solution $\tilde{\xx}$,  $\tilde{\xx}_{\text{opt}}$ is either the optimal solution or the best one known for the instance and $\tilde{\xx}_{\text{init}}$ is the initial solution.

Let $t_{\text{max}} > 0$ be the time limit for executing the heuristic. The primal integral measure is then defined as
\begin{align*}
        P(t_{\text{max}}) = 
        \int_{0}^{t_{\text{max}}} p(t) \, dt.
\end{align*}

\subsubsection{Results}
\label{sec:app1_results}

In order to validate our approach, we compared multiple versions of our \emph{TS-GNN} algorithm with the baselines. Specifically, we compared the following algorithms:
\begin{itemize}
	\item TS with \emph{No-Classifier}, the baseline topology TS algorithm  with an enumeration strategy that evaluates all the possible edge-swap neighbors;
	\item TS with \emph{Random-Add-Classifier}, the baseline topology TS algorithm with a sampling strategy that randomly selects edges to be added for generating edge-swap neighbors;
	\item TS with \emph{Random-Add-Drop-Classifier}, the baseline topology TS algorithm with a sampling strategy that randomly selects both the edges to be dropped and the edges to be added for generating edge-swap neighbors;
	\item TS with \emph{GNN-Add-Classifier}, the TS algorithm plus the GNN classifier for selecting the edges to be added for generating edge-swap neighbors;
	\item TS with \emph{GNN-Drop-Classifier}, the TS algorithm plus the GNN classifier for selecting the edges to be dropped for generating edge-swap neighbors;
	\item TS with \emph{GNN-Add-Drop-Classifier}, the TS algorithm plus the GNN classifiers both for selecting the edges to be dropped and for selecting the edges to be added.
\end{itemize}

\begin{figure}[htbp!]
    \centering
    \includegraphics[width=0.48\linewidth]{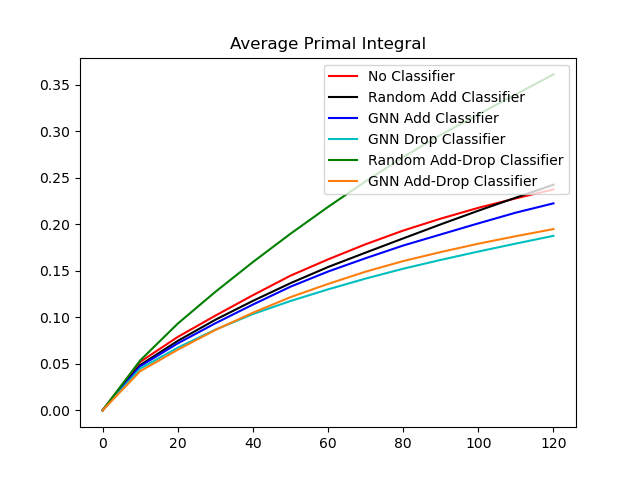}
     \hspace{0.02\linewidth}
    \includegraphics[width=0.48\linewidth]{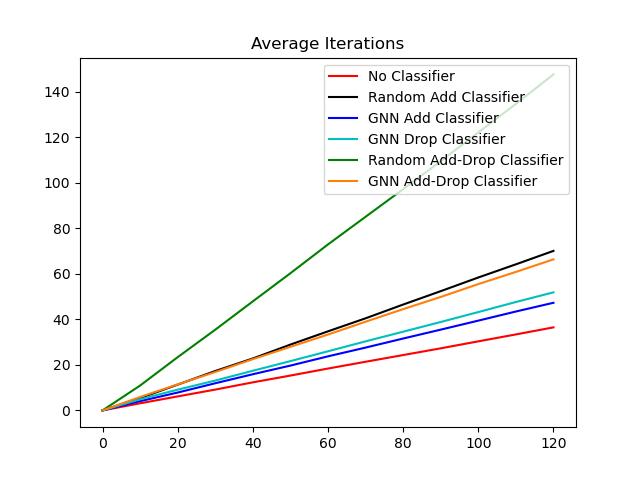}\\
    \caption{Evolution of the average primal integral and the average number of iterations over time on evaluation datasets for the instances of 10 nodes.}
    \label{fig:application1_result1}
\end{figure}

\begin{figure}[htbp!]
    \centering
    \includegraphics[width=0.48\linewidth]{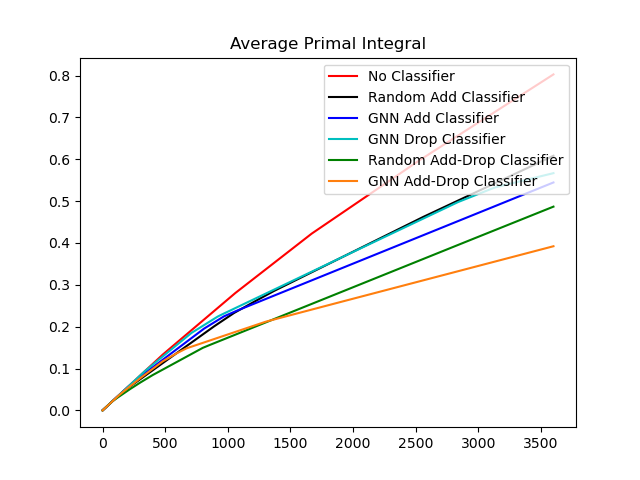}
    \hspace{0.02\linewidth}
    \includegraphics[width=0.48\linewidth]{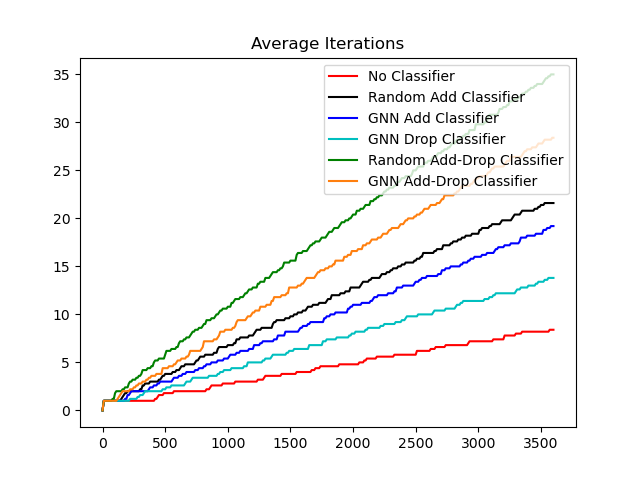}
    \caption{Evolution of the average primal integral and the average number of iterations over time on evaluation datasets for the instances of 30 nodes.}
    \label{fig:application1_result2}
\end{figure}

We evaluated the performance of all the compared algorithms on both small instances with 10 nodes and larger instances with 30 nodes. For each test set, we computed the average primal integral defined in Section \ref{sec:app1_metric} as well as the average number of iterations. The primal integral is our main performance metric for evaluating metaheuristics and it measures the speed of convergence of the objective over the entire search time (the smaller, the better). In addition, the number of iterations counts the number of neighborhoods explored by the algorithm. With the same running time, it reflects how much exploration is guaranteed by modifications made to the basic TS scheme. Indeed, larger number of iterations for the same amount of time indicates a faster moving between neighborhoods and is generally preferable, although the outcome of the search still depends on the quality of the neighborhoods. The results are shown in Figures \ref{fig:application1_result1}  and \ref{fig:application1_result2} for instances with 10 and 30 nodes, respectively. 

From these results, we observe that all our ML-based TS algorithms with GNN classifiers outperform the original TS baseline with no classifier. The \emph{No-Classifier} baseline employs an enumeration strategy that evaluates all the possible moves in the edge-swap neighborhood. Although it is guaranteed to find the best move at each TS iteration, its more extensive exploitation results in the lowest number of iterations over all the compared algorithms. For small instances of 10 nodes, the performance of \emph{No-Classifier} is still decent, as the solution space of these instances is small for exploration. However, it gives the worst results compared to all other algorithms on larger instances because the low number of iterations (slow exploitation) and the resulting lack of exploration becomes the main bottleneck of the algorithm. 

The algorithms with random classifiers, on the contrary, are more efficient in terms of number of iterations. The \emph{Random-Add-Classifier} algorithm explores a random subset of addable edges in the edge-swap neighborhood and the \emph{Random Add-Drop Classifier} generates an even smaller neighborhood by sampling from both droppable and addable edges. On small instances with 10 nodes, \emph{Random-Add-Classifier} achieves a better performance with a lower primal integral although the \emph{Random-Add-Drop-Classifier} algorithm has a larger number of iterations. This is because, for small instances, the quality of the solutions in the sampled neighborhoods is more important than the number of iterations. Whereas for larger instances of 30 nodes, \emph{Random-Add-Drop-Classifier} performs better than \emph{Random-Add-Classifier} since the sampling efficiency becomes more important for exploring larger edge-swap neighborhoods.

On the one hand, the \emph{No-Classifier} baseline can fully exploit each edge-swap neighborhood (\emph{exploitation}) but has a low efficiency in terms of number of iterations. The baselines with random classifiers are more effective in increasing the number of iterations by reducing the size of each neighborhood (\emph{exploration}), however, the quality of the neighborhood is restricted by its sampling efficiency. On the other hand, our complete ML-based algorithm with \emph{GNN-Drop-Add-Classifier} achieves a smaller primal integral than all the baseline algorithms with a reasonably large number of iterations. Moreover, as the size of instances increases, the \emph{GNN-Drop-Add-Classifier} algorithm becomes more competitive and significantly outperforms all the compared algorithms. These results demonstrate that our approach offers a good trade-off in both exploration and exploitation of the solution space. Since the GNN classifiers are only trained with data generated from small instances, the results also show that our method generalizes well on larger instances.

\section{Application 2: Large Neighborhood Search in MIP}
\label{sec:app2}
In this section, we demonstrate how to apply our methodology to another application: a Large Neighborhood Search (LNS) heuristic for solving MIPs.  

\subsection{Large Neighborhood Search}

LNS is a refinement heuristic, i.e., given an initial solution, it is applied to improve the solution by exploring a ``large'' neighborhood. There are several ways of describing a LNS scheme. We adopt the following simple one based of 3 building blocks:
\begin{itemize}
    \item \emph{destroy} function $d$: fixes the values of a subset of variables to the current solution $\boldsymbol{x}$ and ``destroy'' the rest. The output of this function is a sub-MIP with a neighborhood $N(\boldsymbol{x})$. The size of the LNS neighborhood, i.e. the number of variables to be destroyed, will be selected as a hyperparameter;
    \item \emph{repair} function $r$: rebuilds the destroyed solution (in some cases, the repaired solution can be worse than the destroyed solution), typically by solving a sub-MIP defined by $N(\boldsymbol{x})$;
    \item \emph{accept} function $a$: decides whether the new solution should be accepted or rejected. 
\end{itemize}
Given as an input a feasible solution $\boldsymbol{\Bar{x}}$, LNS searches for the best feasible solution $\boldsymbol{\Tilde{x}}$ in the neighborhood of $\boldsymbol{\Bar{x}}$ (the size of the neighborhood is a parameter) and then updates $\boldsymbol{\Bar{x}}\gets \boldsymbol{\Tilde{x}}$. This is repeated until a termination condition is met. The pseudocode of the scheme is given in Algorithm \ref{alg:the_alg_lns}. 
\begin{algorithm}[t]
\caption{The LNS heuristic}
\label{alg:the_alg_lns}
\SetAlgoLined
\KwIn{an initial solution $\xx$;}$\xx^* {\gets} \xx$;

\Repeat {termination condition is reached} 
    {
    $\xx' {\gets} r(d(\xx))$;
    
    \If {$a(\xx', \; \xx)$} { 
         $\xx {\gets} \xx'$;}
    
    \If {$f(\xx) < f(\xx^*)$} { 
         $\xx^* {\gets} \xx$;}
    }

\Return $\xx^*$
\end{algorithm}


\subsection{Learning to Generate Neighbors for LNS}
\label{sec:app2_learn}
 In a LNS heuristic, the solution neighbors are obtained by the destroy function that selects a subset of integer variables to be ``destroyed'', thus freed for neighbor generation. The remaining integer variables will be fixed to their values in the current solution. The classic LNS algorithm applies a randomized sampling strategy or hand-crafted rules for defining the destroy function (see, e.g., \cite{rins}). However, our observation shows that only a small subset of variables have the potential to improve the current solution by changing their values and such a subset is strongly dependent to the structure of the problem. Hence, in order to optimize the performance of LNS, we aim at learning new variable selection strategies to select a subset of variables to be freed. In particular, we investigate the dependencies between the state of the problem, defined by a set of both static and dynamic features collected from the LNS procedure (e.g., context of the problem, incumbent solution) and the binary decisions about integer variables to be either freed or fixed. 
 

We now show how to use a GNN model for learning a destroy policy. The LNS neighborhoods will then be defined based on our pretrained GNN classifier.. 

 Let $\mathcal{S}$ denote the state space of a MIP instance, and let $\mathcal{Y} = \{0,1\}^{|\mathcal{I}|}$ denote all candidate subsets of variables to be destroyed, where $\mathcal{I}$ is the index set of integer variables. We aim to learn a destroy policy $\pi_{\theta^*}$ that maps each state $\ss$ to a label $\yy\in \mathcal{Y}$.  
 
 To generate a training dataset, we consider a set $J$ of instances, and apply $T_j$ iterations of the local branching heuristic \cite{fischetti2003local} to each instance $j\in J$.
We thus determine improving subsets of variables which give us labels. More details about this process can be found in Section \ref{sec: app2_data}.
Let 
$\{\ss_{i}^{(j)}\}_{i=1}^{T_j}$ be the states of the $j^{th}$ instance, and let    $\{\yy_{i}^{(j)}\}_{i=1}^{T_j}$ be the corresponding labels. As explained in Section \ref{sec:learn_train}, the policy $\pi_{\theta^*}$
can be obtained by solving the following classification problem:
\begin{equation}\label{eq:a3_maxlikelihood_3}
    \theta^* = \argmin_{\theta \in \Theta} \sum_{j=1}^{|J|} \sum_{i=1}^{T_j} \mathcal{L} \left(\pi_{\theta}(\ss_{i}^{(j)}), \yy_{i}^{(j)} \right).
\end{equation}

\subsubsection{Feature Design}
We represent each state $\ss$ of a MIP by a bipartite graph $(\mathbf{V}, \mathbf{C}, \mathbf{E})$ \cite{gasse2019exact}.  More precisely, assume that the considered MIP has $n$ variables, $m$ constraints, $d$ features for the variables and $q$ features for the constraints. The variables of the MIP with their feature matrix $\mathbf{V} \in \mathbb{R}^{n\times d}$ are represented on one side of the graph. On the other side are nodes corresponding to the constraints with $\mathbf{C} \in \mathbb{R}^{m\times q}$ being their feature matrix. A variable node $i$ and a constraint node $j$ are connected by an edge $(i,j)$ 
if variable $i$ appears in constraint $j$ of the MIP.
Finally, $\mathbf{E} \in \mathbb{R}^{n\times m\times e}$ denotes the tensor of edge features, with $e$ being the number of features for each edge.
The features in the bipartite graph are listed in Table \ref{table:a3_bipartite-features}.

\begin{table}[htbp!]
    \centering
    \begin{tabular}{c l l}
    \multicolumn{1}{c}{Tensor} & \multicolumn{1}{l}{Feature} & \multicolumn{1}{l}{Description} \\
    \toprule
    \multirow{1}{*}{$\mathbf{V}$}
    & sol\_val & Solution value.  \\
    \midrule
    \multirow{1}{*}{$\mathbf{C}$} & bias & Right-hand side value of the constraint. \\
    \midrule
    \multirow{1}{*}{$\mathbf{E}$}
    & coef & Constraint coefficient. \\
    \bottomrule
    \end{tabular}
    \caption{Description of the features in the bipartite graph $\mathbf{\ss} = (\mathbf{V}, \mathbf{C}, \mathbf{E})$.}
    \vspace{10pt}
    \label{table:a3_bipartite-features}
\end{table}

\subsubsection{GNN Model}
Given that the state of a MIP instance can be represented as a bipartite graph, we propose to use GNN to parameterize the model for the destroy policy. 
Our GNN architecture again consists of 3 modules: the input module, the convolution module, and the output module. For a bipartite graph, a convolution layer is decomposed into two half-layers: one half-layer propagates messages from variable nodes to constraint nodes through edges, and the other one propagates messages from constraint nodes to variable nodes. We refer the reader to \cite{gasse2019exact} for more details.
The output module embeds the features extracted from the convolution module for the prediction of each variable, which maps the graph representation embedding of each variable into a two-neuron output. 

\subsubsection{Loss Function}
The class distribution is highly unbalanced. In the training dataset, we observed that only less than $10\%$ variables belong to the "destroy" class, i.e. the variables in this class have the potential to improve the current solution by changing their values. In order to adapt to the imbalanced distribution, we applied the WCE loss and focal loss \cite{lin2017focal} to train the model.

\subsubsection{LNS Guided by GNNs}
Our refined LNS heuristic with GNN classifier (LNS-GNN) is obtained by using the destroy policy $\pi_{\theta^*}$ in Algorithm \ref{alg:the_alg_lns}.
Hence, given a solution $\xx$ and its associated state $\ss$, we set $\yy$ equal to $\pi_{\theta^*}(\ss)$, we then destroy some variables $\xx_i$ with $\yy_i=1$ to get $N(\xx)$, and we finally obtain a repaired neighbor $\xx'$ of $\xx$ by solving the sub-MIP defined by $N(\xx)$.



\subsection{Numerical Experiments}
\label{sec: app2_experiment}

In this section,  we present the experimental results for the ML-based LNS. As for the WNO application, we first present the data collection, then, we discuss the experimental setting in Section \ref{sec:app2_setting} and the results in Section \ref{sec:app2_results}. The evaluation metrics remain the same as presented in Section \ref{sec:app1_metric}.

\subsubsection{Data Collection}
\label{sec: app2_data}

\paragraph{MIP Benchmark}
To train, evaluate and compare algorithms, we have considered $126$  MIP instances taken from the MIPLIB \cite{gleixner2021miplib} dataset. 
For each instance, an initial feasible solution is required to start the LNS heuristic. We use an intermediate solution found by SCIP \cite{GamrathEtal2020ZR}, typically the best solution obtained by SCIP at the end of the root node computation, i.e., before branching.

\paragraph{Training Data Generation}
To collect data for the classification task, we use the Local Branching (LB) algorithm of Fischetti and Lodi \cite{fischetti2003local}. More precisely, given a MIP instance and an initial incumbent $\bar\xx$, we call the LB algorithm to explore as many neighbors as possible, with a time limit of 600 seconds, and by limiting to $25\%$ the number of destroyed variables.  The best solution $\xx^*$ found by LB is compared with $\bar\xx$ and the variables with changed values are labeled as improving variables for $\xx$. These improving variables are then distroyed from $\xx$, and the partial solution is repaired by calling a MIP solver. The resulting solution becomes the new incumbent to which we repeat the same process. This is done for a set $J$ of instances, as explained in Section \ref{sec:app2_learn}.

\subsubsection{Experimental Setting}
\label{sec:app2_setting}

\paragraph{Training}

In the classification task for producing a good destroy policy, the GNN model learns from the features of the MIP formulation and its incumbent solution. We have trained our GNN classifier with $29$ of the 126 considered  MIP instances.

\paragraph{Evaluation}
The algorithms are compared and evaluated on the remaining $97$ binary MIP instances. 

\paragraph{Experiment Environment}
Our code is written in Python $3.7$ and we use Pytorch $1.60$ \cite{paszke2019pytorch}, Pytorch Geometric $1.7.0$ \cite{fey2019fast}, PySCIPOpt $3.1.1$ \cite{pyscipopt2016}, SCIP $7.01$ \cite{GamrathEtal2020ZR} for developing our models and solving MIPs.

\subsubsection{Results}
\label{sec:app2_results}

We compare our \emph{LNS-GNN} algorithm with the following algorithms:
\begin{itemize}
	\item \emph{LNS-Random}, the baseline LNS algorithm with a sampling strategy that randomly selects variables to destroy;
	\item \emph{LB}, the LB algorithm \cite{fischetti2003local}.
\end{itemize}

For all these algorithms, we have 
limited to 40 the number of variables that can be destroyed to generate neighbors.

\begin{figure}[ht]
          \centering
          \includegraphics[width=0.6\linewidth]{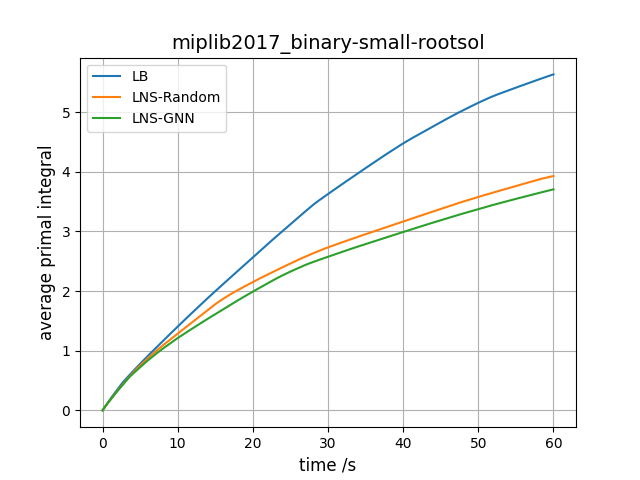}
          \caption{Evaluation results on MIPLIB binary dataset.}
          \label{fig:plot_lns}
\end{figure}

For measuring the heuristic performance of the compared algorithms, we also compute the average primal integral defined in Section \ref{sec:app1_metric}. We ran the listed algorithms for 60 seconds on each instance. The results are shown in Figure \ref{fig:plot_lns}. 


From the results in Figure \ref{fig:plot_lns}, we can see that the primal integral of the \emph{LB} baseline is the largest over the entire solving time, which indicates that exploring all the possible neighbors within the same Hamming distance to the current solution is computationally expensive. Moreover, the \emph{LNS-Random}  baseline performs better than \emph{LB} baseline by applying a random sampling strategy for selecting the subset of variables to be freed. Our \emph{LNS-GNN} algorithm presents the best heuristic behavior in terms of the primal integral, showing that the pretrained GNN model is able to produce structural neighborhoods that contain improving solutions. These results demonstrate that our approach achieves a better trade-off between exploitation and exploration of the promising solution regions. Although, potentially, any of the baseline algorithms could be tuned to obtain better results, this is true for \emph{LNS-GNN} too, and overall we believe that these results show clear promise.

\section{Conclusions and Future Work}
\label{sec:conclusion}

In this work, we presented a methodology for integrating machine learning techniques into metaheuristics for solving combinatorial optimization problems. Namely, we proposed a general ML framework for neighbor generation in metaheuristic search. We firstly defined a neighborhood structure constructed by applying a transformation operator to a selected subset of variables from the current solution. Then, the key of the proposed methodology is to learn a variable selection policy by which the subset transformation operator is able to produce neighborhoods containing high-quality solutions. We formulated the variable selection problem as a classification problem that exploits structural information from the characteristics of the problem and high-quality solutions. 

We demonstrated our methodology on two applications. The first problem we addressed occurs in the context of Wireless Network Optimization, where a Tabu Search metaheuristic is used for the topology design sub-problem. In a predefined topology neighborhood structure, we trained classification models to select sized-reduced, but high-quality neighborhoods. This allowed the metaheuristic search to execute more iterations within the same amount of time. As each iteration requires to solve a series of complex combinatorial sub-problems, more iterations entail a greater exploration of the solution space. In addition, to demonstrate the broader applicability of our approach, we also applied our framework to the Large Neighborhood Search metaheuristic for solving MIPs. The experimental results of the two applications have shown that our approach is able to learn a satisfactory trade-off between the exploration of a larger solution space and the exploitation of promising solution regions in metaheuristic search. 

Although deep neural networks such as GNNs have been wildly applied to represent combinatorial optimization problems, the current GNNs might not be expressive enough to capture all the crucial patterns from data \cite{cappart2021combinatorial}. For future research, it would be interesting to develop more expressive ML models that exhibit better transferability and scalability across broader classes of problems.

\bibliographystyle{unsrt}
\bibliography{references}

\begin{thebibliography}{10}

\bibitem{blum2003metaheuristics}
Christian Blum and Andrea Roli.
\newblock Metaheuristics in combinatorial optimization: Overview and conceptual
  comparison.
\newblock {\em ACM computing surveys (CSUR)}, 35(3):268--308, 2003.

\bibitem{land2010automatic}
Ailsa~H Land and Alison~G Doig.
\newblock An automatic method for solving discrete programming problems.
\newblock In {\em 50 Years of Integer Programming 1958-2008}, pages 105--132.
  Springer, 2010.

\bibitem{Bixby2007}
Robert Bixby and Edward Rothberg.
\newblock Progress in computational mixed integer programming---a look back
  from the other side of the tipping point.
\newblock {\em Annals of Operations Research}, 149(1):37--41, Feb 2007.

\bibitem{bengio2021machine}
Yoshua Bengio, Andrea Lodi, and Antoine Prouvost.
\newblock Machine learning for combinatorial optimization: a methodological
  tour d’horizon.
\newblock {\em European Journal of Operational Research}, 290(2):405--421,
  2021.

\bibitem{talbi2021machine}
El-Ghazali Talbi.
\newblock Machine learning into metaheuristics: A survey and taxonomy.
\newblock {\em ACM Computing Surveys (CSUR)}, 54(6):1--32, 2021.

\bibitem{karimi2022machine}
Maryam Karimi-Mamaghan, Mehrdad Mohammadi, Patrick Meyer, Amir~Mohammad
  Karimi-Mamaghan, and El-Ghazali Talbi.
\newblock Machine learning at the service of meta-heuristics for solving
  combinatorial optimization problems: A state-of-the-art.
\newblock {\em European Journal of Operational Research}, 296(2):393--422,
  2022.

\bibitem{gasse2019exact}
Maxime Gasse, Didier Ch{\'e}telat, Nicola Ferroni, Laurent Charlin, and Andrea
  Lodi.
\newblock Exact combinatorial optimization with graph convolutional neural
  networks.
\newblock In {\em Advances in Neural Information Processing Systems}, pages
  15554--15566, 2019.

\bibitem{gori2005new}
Marco Gori, Gabriele Monfardini, and Franco Scarselli.
\newblock A new model for learning in graph domains.
\newblock In {\em Proceedings. 2005 IEEE International Joint Conference on
  Neural Networks, 2005.}, volume~2, pages 729--734. IEEE, 2005.

\bibitem{scarselli2008graph}
Franco Scarselli, Marco Gori, Ah~Chung Tsoi, Markus Hagenbuchner, and Gabriele
  Monfardini.
\newblock The graph neural network model.
\newblock {\em IEEE transactions on neural networks}, 20(1):61--80, 2008.

\bibitem{hamilton2017representation}
William~L Hamilton, Rex Ying, and Jure Leskovec.
\newblock Representation learning on graphs: Methods and applications.
\newblock {\em arXiv preprint arXiv:1709.05584}, 2017.

\bibitem{cappart2021combinatorial}
Quentin Cappart, Didier Ch{\'e}telat, Elias Khalil, Andrea Lodi, Christopher
  Morris, and Petar Veli{\v{c}}kovi{\'c}.
\newblock Combinatorial optimization and reasoning with graph neural networks.
\newblock {\em arXiv preprint arXiv:2102.09544}, 2021.

\bibitem{dai2017learning}
Elias Khalil, Hanjun Dai, Yuyu Zhang, Bistra Dilkina, and Le~Song.
\newblock Learning combinatorial optimization algorithms over graphs.
\newblock {\em arXiv preprint arXiv:1704.01665}, 2017.

\bibitem{nazari2018reinforcement}
Mohammadreza Nazari, Afshin Oroojlooy, Lawrence Snyder, and Martin Tak{\'a}c.
\newblock Reinforcement learning for solving the vehicle routing problem.
\newblock {\em Advances in neural information processing systems}, 31, 2018.

\bibitem{zhang2020learning}
Cong Zhang, Wen Song, Zhiguang Cao, Jie Zhang, Puay~Siew Tan, and Xu~Chi.
\newblock Learning to dispatch for job shop scheduling via deep reinforcement
  learning.
\newblock {\em Advances in Neural Information Processing Systems},
  33:1621--1632, 2020.

\bibitem{bello2016neural}
Irwan Bello, Hieu Pham, Quoc~V Le, Mohammad Norouzi, and Samy Bengio.
\newblock Neural combinatorial optimization with reinforcement learning.
\newblock {\em arXiv preprint arXiv:1611.09940}, 2016.

\bibitem{gao2020learn}
Lei Gao, Mingxiang Chen, Qichang Chen, Ganzhong Luo, Nuoyi Zhu, and Zhixin Liu.
\newblock Learn to design the heuristics for vehicle routing problem.
\newblock {\em arXiv preprint arXiv:2002.08539}, 2020.

\bibitem{liu2021learning}
Defeng Liu, Andrea Lodi, and Mathieu Tanneau.
\newblock Learning chordal extensions.
\newblock {\em Journal of Global Optimization}, 81(1):3--22, 2021.

\bibitem{he2014learning}
He~He, Hal Daume~III, and Jason~M Eisner.
\newblock Learning to search in branch and bound algorithms.
\newblock {\em Advances in neural information processing systems},
  27:3293--3301, 2014.

\bibitem{Khalil_LeBodic_Song_Nemhauser_Dilkina_2016}
Elias Khalil, Pierre Le~Bodic, Le~Song, George Nemhauser, and Bistra Dilkina.
\newblock Learning to branch in mixed integer programming.
\newblock {\em Proceedings of the AAAI Conference on Artificial Intelligence},
  30(1), Feb. 2016.

\bibitem{khalil2017learning}
Elias Khalil, Hanjun Dai, Yuyu Zhang, Bistra Dilkina, and Le~Song.
\newblock Learning combinatorial optimization algorithms over graphs.
\newblock {\em Advances in neural information processing systems}, 30, 2017.

\bibitem{balcan2018learning}
Maria-Florina Balcan, Travis Dick, Tuomas Sandholm, and Ellen Vitercik.
\newblock Learning to branch.
\newblock In {\em International conference on machine learning}, pages
  344--353. PMLR, 2018.

\bibitem{Liu_Fischetti_Lodi_2022}
Defeng Liu, Matteo Fischetti, and Andrea Lodi.
\newblock Learning to search in local branching.
\newblock {\em Proceedings of the AAAI Conference on Artificial Intelligence},
  36(4):3796--3803, Jun. 2022.

\bibitem{perreault2022}
Vincent Perreault.
\newblock Tactical wireless network design for challenging environments.
\newblock Master's thesis, Ecole Polytechnique, Montreal (Canada), 2022.

\bibitem{paszke2019pytorch}
Adam Paszke, Sam Gross, Francisco Massa, Adam Lerer, James Bradbury, Gregory
  Chanan, Trevor Killeen, Zeming Lin, Natalia Gimelshein, and Luca Antiga.
\newblock Pytorch: An imperative style, high-performance deep learning library.
\newblock {\em Advances in neural information processing systems},
  32:8026--8037, 2019.

\bibitem{fey2019fast}
Matthias Fey and Jan~Eric Lenssen.
\newblock Fast graph representation learning with pytorch geometric.
\newblock {\em arXiv preprint arXiv:1903.02428}, 2019.

\bibitem{berthold2013measuring}
Timo Berthold.
\newblock Measuring the impact of primal heuristics.
\newblock {\em Operations Research Letters}, 41(6):611--614, 2013.

\bibitem{rins}
Emilie Danna, Edward Rothberg, and Claude {Le Pape}.
\newblock Exploring relaxation induced neighborhoods to improve mip solutions.
\newblock {\em Mathematical Programming}, 102:71--90, 2005.

\bibitem{fischetti2003local}
Matteo Fischetti and Andrea Lodi.
\newblock Local branching.
\newblock {\em Mathematical programming}, 98(1-3):23--47, 2003.

\bibitem{lin2017focal}
Tsung-Yi Lin, Priya Goyal, Ross Girshick, Kaiming He, and Piotr Doll{\'a}r.
\newblock Focal loss for dense object detection.
\newblock In {\em Proceedings of the IEEE international conference on computer
  vision}, pages 2980--2988, 2017.

\bibitem{gleixner2021miplib}
Ambros Gleixner, Gregor Hendel, Gerald Gamrath, Tobias Achterberg, Michael
  Bastubbe, Timo Berthold, Philipp Christophel, Kati Jarck, Thorsten Koch, and
  Jeff Linderoth.
\newblock Miplib 2017: data-driven compilation of the 6th mixed-integer
  programming library.
\newblock {\em Mathematical Programming Computation}, pages 1--48, 2021.

\bibitem{GamrathEtal2020ZR}
Gerald Gamrath, Daniel Anderson, Ksenia Bestuzheva, Wei-Kun Chen, Leon Eifler,
  Maxime Gasse, Patrick Gemander, Ambros Gleixner, Leona Gottwald, Katrin
  Halbig, Gregor Hendel, Christopher Hojny, Thorsten Koch, Pierre Le~Bodic,
  Stephen~J. Maher, Frederic Matter, Matthias Miltenberger, Erik M{\"u}hmer,
  Benjamin M{\"u}ller, Marc~E. Pfetsch, Franziska Schl{\"o}sser, Felipe
  Serrano, Yuji Shinano, Christine Tawfik, Stefan Vigerske, Fabian Wegscheider,
  Dieter Weninger, and Jakob Witzig.
\newblock {The SCIP Optimization Suite 7.0}.
\newblock ZIB-Report 20-10, Zuse Institute Berlin, March 2020.

\bibitem{pyscipopt2016}
Stephen Maher, Matthias Miltenberger, Jo{\~a}o~Pedro Pedroso, Daniel Rehfeldt,
  Robert Schwarz, and Felipe Serrano.
\newblock {PySCIPOpt}: Mathematical programming in python with the {SCIP}
  optimization suite.
\newblock In {\em Mathematical Software {\textendash} {ICMS} 2016}, pages
  301--307. Springer International Publishing, 2016.

\end{thebibliography}

\end{document}